\documentclass[12pt]{article}
\usepackage{amsmath, amssymb, amscd, amsthm, amsfonts, xfrac}
\numberwithin{equation}{section}
\usepackage{graphicx}
\usepackage{hyperref}
\usepackage{authblk}
\usepackage{url, cite}
\hypersetup{backref, colorlinks=true}
\newtheorem{Theorem}{Theorem}[section]

\newtheorem{Lemma}[Theorem]{Lemma}

\def \d{\Delta}

\def\qed{\hfill \rule{4pt}{7pt}}
\newenvironment{pf}{{\noindent\bf Proof}\quad}{\hfill $\square$}

\newcommand\Ls[2]{\left( \frac{#1}{#2} \right)}

\date{}
\title{\textbf{Reduction on the congruences of partial sums of P-recursive sequences}}
\author{Qing-Hu Hou and Na Li}
\affil{School of Mathematics \break Tianjin University \break Tianjin 300072, China \break\texttt{qh\_hou@tju.edu.cn, li\_math@tju.edu.cn}}

\begin{document}

\maketitle

\begin{abstract}
Hou and Liu developed a telescoping method to prove the congruence of partial sums of P-recursive sequences. We release the requirement on the telescoper and utilize the congruence of the sequence. With this approach, we are able to confirm a conjecture of Sun and find a new congruence on the central trinomial coefficient.  
\end{abstract}

\noindent {\textit{Keywords:} Motzkin numbers; Central trinomial coefficients; Congruence}


\section{Introduction}
The congruence of partial sums of P-recursive sequences have been widely studied in recent years. These congruences were proved by various methods, especially utilizing combinatorial identities \cite{r1,r2,r3,r4,r5,r6} or symbolic computation \cite{r5,r9, r10}. 

Recently, Hou and Liu \cite{HL8} proposed a telescoping method. Let $a_k$ be a P-recursive sequence of order $d$ and $p(k), q(k)$ be two polynomials with integral coefficients and $q(0)=0$.
To show the congruence
\begin{equation}\label{eq-pak}
\sum_{k=0}^{n-1} p(k) a_k \equiv 0 \pmod{q(n)},
\end{equation}
they try to find polynomials with integral coefficients $f_0(k), \ldots, f_{d-1}(k)$ such that
\[
p(k) a_k = \Delta_k \left( q(k) \sum_{i=0}^{d-1} f_i(k) a_{k-i} \right),
\]
where $\Delta_k f(k) = f(k+1)-f(k)$ is the difference operator.

In this paper, we further develop this method. We release the requirement that $q(k)$ appears as a factor in the difference and try to find polynomials with integral coefficients $g_0(k), \ldots,g_{d-1}(k)$ such that
\[
p(k) a_k = \Delta_k \left( \sum_{i=0}^{d-1} g_i(k) a_{k-i} \right).
\]
Then we will show 
\[
 \sum_{i=0}^{d-1} g_i(n) a_{n-i} - \sum_{i=0}^{d-1} g_i(0) a_{-i} \equiv 0 \pmod{q(n)},
\]
which clearly leads to \eqref{eq-pak}. To this end, we reduce $g_i(n)$ to $\tilde{g}_i(n)$ with the property that
\[
g_i(n) \equiv \tilde{g}_i(n) \pmod{q(n)}
\]
and utilize the divisibility of $a_n$.

We mainly consider the product $c_n$ of two P-recursive sequences $a_n$ and $b_n$. Although $c_n$ is also P-recursive, but it is more flexible to use $a_n$ and $b_n$. For example, suppose that the order of $a_n$ and $b_n$ are $d,e$ respectively. We will seek for polynomials $p(k)$ and $g_{i,j}(k)$ such that
\[
p(k) a_k b_k  = \Delta_k \left( \sum_{i=0}^{d-1} \sum_{j=0}^{e-1} g_{i,j}(k) a_{k-i} b_{k-j} \right).
\]

We have implemented a {\tt Maple} package {\tt RC} to accomplish the computations. We will illustrate the method by three examples. 

As the first example, we confirm the conjecture of Sun \cite[Conjecture 5.1(i)]{r1}.
\begin{Theorem}\label{thW}
Let 
\begin{equation}\label{eq-W}
	W_n = \sum\limits_{k=0}^{\lfloor n/2 \rfloor}\binom{n}{2k}\frac{\binom{2k}{k}}{2k-1}.
\end{equation}

	\begin{enumerate}
		\item[\rm(i)] For any $n\in\mathbb{Z}^+$ ,we have
		\begin{align}\label{W_kn}
			\sum\limits_{k=0}^{n-1}(8k+9)W_k^2\equiv n \ ({\rm mod} \ 2n).
		\end{align}\par
		\item[\rm(ii)] For any odd prime $p$, we have 
		\begin{align}\label{W_kp2}
			\frac{1}{p}\sum\limits_{k=0}^{p-1} (8k+9)W_k^2 \equiv 24+10 \Ls{-1}{p} -9 
		\Ls{p}{3} - 18 \Ls{3}{p} \pmod{p},
		\end{align}
	where $\Ls{\cdot}{\cdot}$ is the Legendre symbol.
	\end{enumerate}
\end{Theorem}

The other two examples are related to the $n$-th trinomial coefficient $T_n$, which is given by
\begin{align*}
	T_n = [x^n] (1+x+x^2)^n = \sum_{l=0}^{\lfloor k/2 \rfloor}\binom{k}{2l}\binom{2l}{l}=\sum_{l=0}^{k}\binom{k}{l}\binom{k-l}{l}.
\end{align*}

By use of combinatorial identities, Sun \cite{r1} showed that
\begin{equation} \label{con-Sun}
\frac{n^2(n^2-1)}{6} \ \left| \ \sum_{k=0}^{n-1}k(k+1)(8k+9)T_kT_{k+1}. \right.
\end{equation}
We give a new proof and further derive 
\begin{Theorem} \label{pTkTk+1}
For any prime $p>3$, we have
\begin{equation} \label{eq-pTk+1}
\sum_{k=0}^{p-1}k(k+1)(8k+9)T_kT_{k+1}\equiv -p^2 \left( \frac{53}{12} + \frac{21}{4} \Ls{p}{3} \right) \pmod{p^3}.
\end{equation}
\end{Theorem}

Finally, we discover the following congruences.
\begin{Theorem}\label{TkTk-1}
	\begin{enumerate}
		\item[{\rm (i)}] For any $n\in \mathbb{Z}^+$, we have
		\begin{equation}\label{nTkTk-1}
		2\sum_{k=0}^{n-1} (k+1)(16k+21)T_kT_{k+1} \equiv 9 nT_{n-1}T_n\ ({\rm mod} \ n^2).
		\end{equation}
	
		\item[{\rm (ii)}] For any prime $p\neq 3$, we have 
		\begin{multline} \label{p2TkTk-1}
		\sum_{k=0}^{p-1} (k+1)(16k+21)T_kT_{k+1} \\
		\equiv p \frac{3^{p+1}}{2} \Ls{p}{3} + p^2 \left( \frac{17}{4} + \frac{57}{4} \Ls{p}{3} \right)  \pmod{p^3}.
		\end{multline}
	\end{enumerate}
\end{Theorem}

\section{Proof of  Theorem \ref{thW}}

Our package {\tt RC} is based on the linear recurrence of the sequence. From the definition \eqref{eq-W} of $W_n$, we can easily derive its recurrence relation by Zeilberger's algorithm. In fact, we have
\begin{align}\label{recurW}
(n+3)W_{n+3}=-3(n+1)W_{n}+(n-5)W_{n+1}+(3n+7)W_{n+2}.
\end{align}
To make  the recurrence relation \eqref{recurW} holds for any integer $n$, we extend the definition of $W_n$ to $n<0$ by setting
\[
W_{n} = \begin{cases}
-1, & \mbox{if $n=-1$}, \\
0, & \mbox{otherwise}.
\end{cases}
\]
Now using the package {\tt RC}, we find that 
\begin{align*}
(8n+9)W_n^2=\d_n S_n,
\end{align*}
with 
\begin{multline*}
S_n =\frac{9}{2}n(n-1)W_{n-1}W_{n-2}-\frac{3}{2}n(n-1)W_{n-2}W_{n} +\frac{1}{2}n(n-1)W_{n-1}W_{n} \\
- \frac{3}{2}n(n-7)W_{n-1}^2.
\end{multline*}
We thus derive
\begin{equation}\label{eq-WS}
\sum\limits_{k=0}^{n-1}(8k+9)W_k^2 = S_n.	
\end{equation}

With \eqref{eq-WS}, we can prove the first part of Theorem~\ref{thW}. We first reduce $S_n$ modulo $2n$.
\begin{align*}
	S_n & \equiv \frac{1}{2}n(n-1)W_{n-1}W_{n-2}+\frac{1}{2}n(n-1)W_{n-2}W_{n}+\frac{1}{2}n(n-1)W_{n-1}W_n\\
	&\quad +\frac{1}{2}n(n+1)W_{n-1}^2 \pmod{2n}\\
	&=\frac{1}{2}n(n-1)(W_{n-1}+W_{n})(W_{n-2}+W_{n-1})+nW_{n-1}^2.
\end{align*}
Noting that ${2k \choose k} = 2 {2k-1 \choose k-1}$, we see that $W_n$ is odd for $n \ge 0$. Therefore, 
\[
	nW_{n-1}^2 \equiv  n \pmod{2n}
\]
and 
\[
(W_{n-1}+W_n) (W_{n-1}+W_{n-2}) \equiv 0 \pmod{4}, \quad n \ge 2,
\]	
implying \eqref{W_kn}.

To prove \eqref{W_kp2}, we need evaluate $W_p,W_{p-1}$ and $W_{p-2}$ modulo $p$.
\begin{Lemma}\label{lemmWp}
Let $p>3$ be a prime. Denote $a=\Ls{-3}{p}$ and $b=\Ls{-1}{p}$. We have
\begin{align}
W_{p} & \equiv -1 - p (1+ 3a - 4b) \pmod{p^2}, \\
W_{p-1} & \equiv 3a-4b \pmod{p}, \\
W_{p-2} & \equiv 7a-8b \pmod{p}.
\end{align}
\end{Lemma}

\pf	
All of these three congruences are related to 
\[
\sum_{k=0}^{\frac{p-1}{2}}  {2k \choose k} \pmod{p} \quad \mbox{and} \quad \sum_{k=0}^{\frac{p-1}{2}} \frac{1}{2k-1} {2k \choose k} \pmod{p}.
\]
So we evaluate them first. 

It is well-known (see, for example, Theorem 1.2 and Lemma 3.1 in \cite{PS}) that 
\[
\sum_{k=0}^{\frac{p-1}{2}} {2k \choose k} \equiv \Ls{p}{3} \equiv a \pmod{p}
\]
and
\[
{p-1 \choose \frac{p-1}{2}} \equiv (-1)^{\frac{p-1}{2}}  = b \pmod{p}.
\]
By Extended Zeilberger’s algorithm { \cite{EZ}}, we fnd that
\begin{align*}
	3\binom{2k}{k}-\frac{1}{2k-1}\binom{2k}{k}=\Delta_k \left( \frac{2k}{2k-1}\binom{2k}{k} \right). 
\end{align*}
So
\begin{align*}
\sum_{k=0}^{\frac{p-1}{2}} \frac{1}{2k-1}\binom{2k}{k} & =  3 \sum_{k=0}^{\frac{p-1}{2}} {2k \choose k} - \frac{p+1}{p}  {p+1 \choose \frac{p+1}{2}} \\
& =3 \sum_{k=0}^{\frac{p-1}{2}} {2k \choose k} - 4  {p-1 \choose \frac{p-1}{2}} 
\equiv 3a - 4 b \pmod{p}.
\end{align*}

By Theorem 1.2 of \cite{PS}, we have
\begin{equation}\label{eq-bi/k}
\sum_{k=1}^{\frac{p-1}{2}} \frac{{2k \choose k}}{k} \equiv \sum_{k=1}^{p-1} \frac{{2k \choose k}}{k} \equiv 0 \pmod{p}.
\end{equation}
Hence
\begin{align*}
	W_p 	& = - 1 + \sum_{k=1}^{\frac{p-1}{2}} \frac{1}{2k-1}\binom{p}{2k}\binom{2k}{k} \\
	&\equiv -1 -\sum_{k=1}^{\frac{p-1}{2}}\frac{p}{2k(2k-1)} \binom{2k}{k}  \\
	&= -1 + p \sum_{k=1}^{\frac{p-1}{2}} \frac{\binom{2k}{k}}{2k}-p\sum_{k=1}^{\frac{p-1}{2}}\frac{1}{2k-1}\binom{2k}{k}\\
	&\equiv -1 - p (1+3a-4b) \pmod{p^2}.
\end{align*}
We also have
\begin{align*}
	W_{p-1}&\equiv \sum_{k=0}^{\frac{p-1}{2}}\binom{2k}{k}\frac{1}{2k-1} = 3a-4b \pmod{p}.
\end{align*}	
Finally, we have
\begin{align*}
W_{p-2} & = \sum_{k=0}^{\frac{p-3}{2}}  \frac{1}{2k-1} \frac{p-1-2k}{p-1} {p-1\choose 2k} {2k \choose k} \\
& \equiv  \sum_{k=0}^{\frac{p-3}{2}}  \frac{2k+1}{2k-1} {2k \choose k} \\
& = \sum_{k=0}^{\frac{p-1}{2}} {2k \choose k} - {p-1 \choose \frac{p-1}{2}}
	+ 2 \sum_{k=0}^{\frac{p-1}{2}} \frac{1}{2k-1} {2k \choose k} - \frac{2}{p-2} {p-1 \choose \frac{p-1}{2}} \\
& \equiv 7a - 8b \pmod{p},
\end{align*}
completing the proof of Lemma~\ref{lemmWp}. \qed

Now we are ready to prove the second part of Theorem~\ref{thW}. When $p=3$, we can directly verify that (\ref{W_kp2}) holds. So we assume that $p>3$ in the following context. Setting $n = p-1$ in \eqref{eq-WS}, we derive that
\[
\frac{1}{p}	\sum_{k=0}^{p-1}(8k+9)W_k^2
\equiv -\frac{9}{2}W_{p-1}W_{p-2}+\frac{3}{2}W_{p-2}W_{p}-\frac{1}{2} W_{p-1}W_p
+ \frac{21}{2}W_{p-1}^2 \pmod{p}.
\]
Substituting the evaluations of $W_p, W_{p-1}$ and $W_{p-2}$, we deduce that
\[
\frac{1}{p}	\sum_{k=0}^{p-1}(8k+9)W_k^2 \equiv 24b^2 -18ab-9a+10b \pmod{p},
\]
implying \eqref{W_kp2}.

\section{Proofs of Theorems~\ref{pTkTk+1} and \ref{TkTk-1}}

By our package {\tt RC}, we find that
\[
k (k+1) (8k+9) T_k T_{k+1} = \Delta S_k,
\]
with
\begin{align}
S_k & = - \frac{1}{24}  k^2  (2k-5)^2  T_k^2 + \frac{1}{4} k^2 (4 k^2 + 20 k - 21)  T_{k-1}T_k \nonumber \\
& \qquad - \frac{3}{8} k^2 (2k-3)^2 T_{k-1}^2. \label{Sk1} \\
& = -\frac{1}{24}(k+1)^2(2k-3)^2T_{k+1}^2+\frac{1}{4}(k+1)(4k^3-5k+3)T_{k}T_{k+1} \nonumber \\
& \qquad -\frac{3}{8}(k+1)^2(2k-1)^2 T_{k}^2 \label{Sk2}.
\end{align}

Let
\[
A_n = 6 \sum_{k=0}^{n-1} k (k+1) (8k+9) T_k T_{k+1}.
\]
By \eqref{Sk1}, we have
\begin{align*}
A_n & =  - \frac{1}{4}  n^2  (2n-5)^2  T_n^2 + \frac{3}{2} n^2 (4 n^2 + 20 n - 21)  T_{n-1}T_n  - \frac{9}{4} n^2 (2n-3)^2 T_{n-1}^2 \\ 
&  \equiv - n^2 \left( \frac{T_n + T_{n-1}}{2} \right)^2 \pmod{n^2}.
\end{align*}
Since 
\[
T_k = 1 + 2 \sum_{\ell=1}^{\lfloor k/2 \rfloor} {k \choose 2\ell} {2 \ell -1 \choose \ell-1}
\]
is odd, we immediately derive that $n^2 \mid A_n$. Similarly, by \eqref{Sk2} we deduce that
\begin{align*}
A_n & \equiv -\frac{1}{4} (n+1)^2 T_{n+1}^2 + \frac{3}{2} (n+1)^2 T_{n}T_{n+1} - \frac{9}{4} (n+1)^2 T_{n}^2 \pmod{n^2 - 1} \\
&= - (n^2-1) \left(\frac{T_{n+1}-3T_n}{2} \right)^2 - (n+1) \frac{(T_{n+1}-3T_n)^2}{2}.
\end{align*}
Since $T_n$ are odd, the first part is divisible by $n^2-1$. For the second part, we will show that 
\begin{equation}\label{dif3}
T_{n+1}- 3T_n \equiv 0 \pmod{n-1}.
\end{equation}
In fact, let $t_n = \frac{T_{n+1}-3T_n}{2}$. By Zeilberger's algorithm, we have
\[
3(n-1) t_{n-2} +2(n-1) t_{n-1} - (n+1) t_n = 0.
\]
Modulo $n-1$, we derive that
\[
2 t_n \equiv 0 \pmod{n-1}.
\]
Therefore \eqref{dif3} holds and thus $(n^2-1) \mid A_n$. Combining with the fact that $n^2 \mid A_n$, we reprove the congruence \eqref{con-Sun}.

Now we give the proof of Theorem~\ref{pTkTk+1}. 

\pf
We start from the expression \eqref{Sk1} of $S_k$. Let $n=p$ be a prime number greater than $3$. We have
\[
\frac{1}{p^2} \sum_{k=0}^{p-1} k (k+1) (8k+9) T_k T_{k+1}
 \equiv -\frac{25}{24} T_p^2 -\frac{21}{4} T_{p-1} T_p - \frac{27}{8} T_{p-1} \pmod{p}. 
\]
It has been shown by Sun \cite[Lemma 2.5]{r2} that  
\[
T_p \equiv 1 \pmod{p} \quad \mbox{and} \quad T_{p-1} \equiv \Ls{-3}{p} = \Ls{p}{3} \pmod{p}.
\]
We immediately derive \eqref{eq-pTk+1}. \qed

At last, we prove Theorem~\ref{TkTk-1}.

\pf Once again, by the package {\tt RC} we find that
\begin{multline}
S_n = \sum_{k=0}^{n-1} (k+1)(16k+21) T_k T_{k+1} = - \frac{1}{8} n^2(4n-7) T_n^2 + \frac{3}{4} n (4n^2+19n+6) T_{n-1}T_n \\
 - \frac{9}{8} n^2(4n-3) T_{n-1}^2. \label{eq-d2}
\end{multline}
Then , we get
\begin{align*}
2S_n & \equiv -\frac{1}{4} n^2T_n^2 + \frac{1}{2} n (n+18) T_{n-1} T_n - \frac{1}{4} T_{n-1} \\
& = -n^2 \left( \frac{T_n - T_{n-1}}{2} \right)^2 + 9n T_{n-1}T_n.
\end{align*}
Noting that the first part is divisible by $n^2$, we immediately get \eqref{nTkTk-1}.

When $n=p$ is a prime number greater than $3$, we have
\[
S_n  \equiv \frac{7}{8} p^2 T_p^2 + \frac{3}{4} p (19p+6) T_{p-1}T_p 
+ \frac{27}{8} p^2 T_{p-1}^2 \pmod{p^3}.
\]
By \eqref{eq-bi/k}, we have
\begin{equation}\label{eq-Tpm2}
T_p-1= \frac{p}{2} \sum_{k=1}^{\frac{p-1}{2}} \frac{1}{k} \binom{p-1}{2k-1}\binom{2k}{k}
\equiv - \frac{p}{2}\sum_{k=1}^{\frac{p-1}{2}} \frac{\binom{2k}{k}}{k}
 	\equiv 0 \pmod{p^2}.
\end{equation}
In \cite{r7}, Sun showed that for any $m\in \mathbb{Z}$ and $n\in \mathbb{Z}^+$, 
 \begin{equation}\label{eq-bitr}
 	\sum_{k=0}^{n-1}\binom{n-1}{k}\binom{2k}{k}(-1)^km^{n-1-k} = \sum_{k=0}^{\lfloor\frac{n-1}{2} \rfloor}\binom{n-1}{k}\binom{n-1-k}{k}(m-2)^{n-1-2k}
 \end{equation}
 and 
 \begin{align*}
 	\sum_{k=0}^{p-1}\binom{p-1}{k}\frac{\binom{2k}{k}}{(-3)^k}\equiv \Ls{p}{3} \pmod{p^2}.
 \end{align*}
Setting $m=3$ and $n=p$ in \eqref{eq-bitr}, we derive that
\begin{align}
T_{p-1}  = \sum_{k=0}^{\frac{p-1}{2}} {p-1 \choose k} {p-1-k \choose k} 
& = 3^{p-1}\sum_{k=0}^{p-1}\binom{p-1}{k}\frac{\binom{2k}{k}}{(-3)^k} \nonumber \\
& \equiv 3^{p-1} \Ls{p}{3} \pmod{p^2}. \label{eq-Tp-1m2}
\end{align}
Substituting \eqref{eq-Tpm2} and \eqref{eq-Tp-1m2}, we immediately derive \eqref{p2TkTk-1}. \qed



\begin{thebibliography}{abc99xyz}
	

	\bibitem{PS}
	H. Pan and Z.W. Sun, A combinatorial identity with application to Catalan numbers, Discrete Math. {\bf 306(16)} (2006) 1921--1940.
	\bibitem{r4}H.-Q Cao, H. Pan. Some congruences for trinomial coefficients, Houston J  Math 40 (2014), 1073-1087.
	\bibitem{r3}J.-C. Liu, Supercongruences involving Motzkin numbers and central trinomial coefficients. ArXiv: 2208.10275 (2022).
	\bibitem{HL8}
	Q.-H. Hou and K. Liu, Congruences and telescopings of P-recursive sequences, J. Diff. Equ. Appl. {\bf 27(5)} (2021) 686--697.
	\bibitem{HL} Q.-H. Hou and K. Liu, Congruences and telescopings of P-recursive sequences, J. Differ. Equ. Appl. (2021) 1--12.
	\bibitem{r9} Q.-H. Hou, Y.-P. Mu and D. Zeilberger. Polynomial reduction and supercongruences. J. Symbolic Comput., 103(2021), 127–140.
	\bibitem{r10} R.-H. Wang, M X X. Zhong. Polynomial Reduction for Holonomic Sequences. arXiv: 2205.11129( 2022).
	\bibitem{EZ} William Y.C. Chen., Q.-H. Hou, Y.-P. Mu, 2012. The extended Zeilberger algorithm with parameters. J. Symb. Comput. 47 (6), 643–654.
	\bibitem{r5}Y.-P. Mu, Z.-W. Sun, Telescoping method and congruences for double sums, Int. J. Number Theory	14 (2018) 143–165.
	\bibitem{r6}Z.-W. Sun, ARITHMETIC PROPERTIES OF APÉRY NUMBERS AND CENTRAL DELANNOY NUMBERS, Arxiv: 1006.2776 (2010).
	\bibitem{r2}Z.-W. Sun, Congruences involving generalized central trinomial coefficients, Sci. China	Math. 57 (2014), 1375–1400.
	\bibitem{r7}Z.-W. Sun, On sums of binomial coefficients modulo p2, Colloq. Math. 127 (2012) 39–54
	\bibitem{r1}Z.-W. Sun, On Motzkin numbers and central trinomial coefficients, Adv. in Appl. Math. 136 (2022), Art. 102319.


	
\end{thebibliography}
\end{document}